\documentclass[12pt]{article}
\usepackage{amssymb}
\usepackage{latexsym}
\usepackage{amsmath}
\usepackage{amscd}
\newfont{\cyr}{wncyr10}
\newfont{\cyb}{wncyb10}

\begin{document}

\setlength{\baselineskip}{24pt}
\setlength{\parindent}{.5in}
\setlength{\parskip}{0pt}
\setlength{\footskip}{.3in}

\setlength{\textwidth}{6.3in}
\setlength{\oddsidemargin}{0in}
\setlength{\topmargin}{-.5in}
\setlength{\textheight}{9in}

\def\cc {{\mathfrak c}}
\def\ii {{\mathfrak i}}
\def\UU {{\mathfrak U}}
\def\CC {{\Bbb C}}
\def\HH {{\Bbb H}}
\def\NN {{\Bbb N}}
\def\PP {{\Bbb P}}
\def\QQ {{\Bbb Q}}
\def\RR {{\Bbb R}}
\def\TT {{\Bbb T}}
\def\ZZ {{\Bbb Z}}
\def\sA {{\mathcal A}}
\def\sB {{\mathcal B}}
\def\sC {{\mathcal C}}
\def\sD {{\mathcal D}}
\def\sE {{\mathcal E}}
\def\sF {{\mathcal F}}
\def\sG {{\mathcal G}}
\def\sH {{\mathcal H}}
\def\sI {{\mathcal I}}
\def\sJ {{\mathcal J}}
\def\sK {{\mathcal K}}
\def\sL {{\mathcal L}}
\def\sM {{\mathcal M}}
\def\sN {{\mathcal N}}
\def\sO {{\mathcal O}}
\def\sP {{\mathcal P}}
\def\sQ {{\mathcal Q}}
\def\sR {{\mathcal R}}
\def\sS {{\mathcal S}}
\def\sT {{\mathcal T}}
\def\sU {{\mathcal U}}
\def\sV {{\mathcal V}}
\def\sW {{\mathcal W}}
\def\sX {{\mathcal X}}
\def\sY {{\mathcal Y}}
\def\sZ {{\mathcal Z}}
\def\od {\mathrm{od}}
\def\cf {\mathrm{cf}}
\def\dom {\mathrm{dom}}
\def\id {\mathrm{id}}
\def\int {\mathrm{int}}
\def\cl {\mathrm{cl}}
\def\Hom {\mathrm{Hom}}
\def\ker {\mathrm{ker}}
\def\log {\mathrm{log}}
\def\nwd {\mathrm{nwd}}

\def\TT{{\Bbb T}}
\def\Z{{\Bbb Z}}
\def\ZZ{{\Bbb Z}}

\newtheorem{thm}{Theorem}[section]

\newtheorem{theorem}[thm]{Theorem}
\newtheorem{corollary}[thm]{Corollary}
\newtheorem{lemma}[thm]{Lemma}
\newtheorem{claim}[thm]{Claim}
\newtheorem{axiom}[thm]{Axiom}
\newtheorem{conjecture}[thm]{Conjecture}
\newtheorem{fact}[thm]{Fact}
\newtheorem{hypothesis}[thm]{Hypothesis}
\newtheorem{assumption}[thm]{Assumption}
\newtheorem{proposition}[thm]{Proposition}
\newtheorem{criterion}[thm]{Criterion}
\newtheorem{definition}[thm]{Definition}
\newtheorem{definitions}[thm]{Definitions}
\newtheorem{discussion}[thm]{Discussion}
\newtheorem{example}[thm]{Example}
\newtheorem{notation}[thm]{Notation}
\newtheorem{remark}[thm]{Remark}
\newtheorem{remarks}[thm]{Remarks}
\newtheorem{problem}[thm]{Problem}
\newtheorem{principle}[thm]{Principle}
\newtheorem{question}[thm]{Question}
\newtheorem{questions}[thm]{Questions}
\newtheorem{notation-definition}[thm]{Notation and Definitions}
\newtheorem{acknowledgement}[thm]{Acknowledgement}

\title{
Introduction to typed topological spaces 
\footnote{
%, 06B05, 62A01
;  
Secondary  
%54J99, 
06B99, 62P99}
\\
To the memory of W.W. Comfort 
%\footnote{Thanks}
%\subjclass{Primary 54A04, 06B05, 62A01;  Secondary  54J99, 06B99, 62P99}
%\keywords{Partital Order Set, Topology
}

\author{Wanjun Hu
\\Department of Math \& CS\\ Albany State University\\ Albany, GA 31705}
\date{November 3, 2013}
\maketitle

\newcommand{\s }{\mathcal }
\newcommand{\bA}{\mathbf{A}}
\newcommand{\bB}{\mathbf{B}}
\newcommand{\bC}{\mathbf{C}}
\newcommand{\bD}{\mathbf{D}}
\newcommand{\bI}{\mathbf{I}}
\newcommand{\bE}{\mathbf{E}}
\newcommand{\bK}{\mathbf{K}}
\newcommand{\bT}{\mathbf{T}}

\begin{abstract}
The concept of typed topological space is introduced, for which open sets in a topology on a finite set will be assigned types (from lattice). The neighborhood system of a point, the closure and the connectedness can be defined according to chain of types, which effectively avoids the situation when most singletons are closed an open. Furthermore, statistics can be used to provide semantics of points with statistical characteristics. 


\end{abstract}

\maketitle

\section{Introduction}

\indent\indent Topological spaces in general study properties of open sets under arbitrary union and finite intersection. The concept of neighborhood, separation axioms and covering properties help describe the local and global behaviors of a space. Spaces with separation axioms and covering properties are useful tools to study infinite spaces. A finite space with $T_1$ separation property is automatically discrete, which renders such kind of spaces uninteresting. 

Finite topological spaces have been discussed in many areas such as Psychology (\cite{lewin}), Social Science, Database (and big data). Because of the nature of being discrete, there has not had any substantial applications of tools and concepts from General Topology in those areas. However, the concepts of neighborhood, approximation, and limit are still widely used.  

Similarly, a database contains several tables, and each table contains a fixed number of columns and arbitrary length of rows. When treating each row as a point and each column as a predicate to define open sets, which will have a topological space. The SQL commands let us to retrieve information by selecting data using union and intersections. The space can easily be discrete. 

On the other hand, one may have several different ways to define open sets on a given finite set. For example, in a community with several streets. We can use left-neighbor and right-neighbor on the same street to define open sets, which eventually becomes a discrete topology. We can also use use friendship, relative and other properties to define open sets, which will becomes a discrete topology too. The topology are the same, but  there do have something different here.

To facilitate applications of general topological concepts and tools to above mentioned examples, 
we define a so-called typed topology on finite sets, in which each open set is associated 
with a type. Set inclusion of open sets will assume orders on types. Neighborhood systems of points and closure operations will be restricted to chains of types. Furthermore, statistical measures can be used to provide semantics to points and pairs of points.

In Section 2, we define the concept of typed topological spaces. In Section 3, we study neighborhood systems of points restricted to chains of types. In Section 4, we study the closure operations and density property restricted to chains of types. In Section 5, we introduce connectedness under chains of types and statistical measures.

  Concepts from general topology follow from (\cite{englking}). 
%    Given a set $X$, a topology $\sT$ (\cite{englking}) is a family of subsets of $X$
%satisfying: (1) $\emptyset\in \sT$; (2) $X\in\sT$; (3) for any $U_1, U_2\in\sT$, $U_1\cap U_2\in\sT$; %and (4) for any subfamily $\sT'\subseteq\sT$, the union $\bigcup\{U: U\in\sT'\}\in\sT$.
%The pair $(X, \sT)$ is called a $topological$ $space$. A topology can easily be obtained
%from a family $\sB$ of subsets of $X$ by the operations of finite intersection and arbitrary
%unions of elements in $\sB$, in which case $\sB$ is called a $subbase$.
         For a set $P$, the binary relation $\leq$ is called  a $partial$ $order$ (\cite{birkhoff}) 
if for any elements $a, b, c\in P$, the following conditions hold: (1) $a\leq a$; (2) $a\leq b$ and $b\leq a$ imply $a=b$; and (3) $a\leq b$ and $b\leq c$ imply $a\leq c$. The notation $a<b$ is used when $a\leq b$ and $a\neq b$. The pair $(P, \leq)$ is a lattice if for any two elements $a,b\in P$, both the greatest lower bound (the meet) $a\wedge b$ and least upper bound (the join) $a\vee b$ exist. A lattice is bounded if bottom $0$ and top $1$ exist. Furthermore, a lattice is distributive if the two operators $\wedge$ and $\vee$ obey the distribution law.
        Two elements $p, q\in P$ 
are called compatible if there exists an element $r\in P$ satisfying $r\leq p$ and $r\leq q$. 
A proper subset $F\subset P$
 is called a $filter$ if (1) $0\notin F$ provided $0$ exists; (2) for any $p,q\in F$, there exists $r\in F$ satisfying 
$r\leq p$ and $r\leq q$;
 and (3) for any $p,q\in P$ with $p\leq q$, $q\in F$ whenever $p\in f$.

In this paper, we assume all topological spaces, partially ordered sets are finite and all lattices are  bounded finite distributive lattices. 

\section{Typed Topology}

\indent\indent Let $X$ be a set and let $(P, \leq)$ be a partially ordered set.  For each $x\in X$, we define two element, $x's$ and $\neg x's$, which means belonging to $x$ and not belonging to $x$ respectively. We denote the set $\{x's: x\in X\}$ by $X's$ and $\{\neg x's: x\in X\}$ by $\neg X's$. We consider $X's\cup\neg X's$ as be unordered. Furthermore, The set $P\cup X's\cup \neg X's$ is denoted $L_0$. We consider the partially ordered set $(L_0, \leq)$ as an extension of $(P, \leq)$ in the way that $\leq$ restricted to $P$ is $(P, \leq)$, and $\leq$ restricted to $X's\cup\neg X's$ is only $x\leq x$ for any element $x$.

\indent Following (\cite{dean}), the free distributive lattice  $L(P, X)$ generated by $L_0$ is defined as the family of terms (or words) in $L_0$ in the form of $A_1\vee A_2\vee ... \vee A_k$, where each $A_i$ ($1<i\leq k$) 
is a meet of elements in $L_0$. Furthermore, $(L(P,X)$ has a partial order that extends $(P, \leq)$, which is defined as follows (\cite{dean2}). The order relation $A\leq B$ holds for $A, B\in L(P, X)$ when one of the following is true:  (1)  for some $a,b\in L_0$, $A=a, B=b$ and $a\leq b$; (2) for some $S\subseteq L_0$ and $b\in L_0$ with $A=\wedge\{a:a\in S_2\}$ and $B=b$, $a\leq b$ holds for some $b\in S$; (3) for some $S_1, S_2\subseteq L_0$ with $A=\wedge\{a: a\in S_1\}$ and $B=\wedge\{b: b\in S_2\}$, $A\leq b$ holds for all $b\in S_2$; (4)for $A=A_1\vee A_2\vee ... \ A_m$ and $B=\wedge\{b:b\in S\}$ for some $S\subseteq L_0$, where $A_i's$ are meets of elements in $L_0$, $A_i\leq B$ holds for all $1\leq i\leq m$; (5) for $B=B_1\vee B_2\vee ... \vee B_n$, where $B_j's$ are meets of elements in $L_0$, $A\leq B_j$ holds for all $1\leq j\leq n$. 
    Furthermore, we consider $L(P, X)$ as being bounded in the way that  for any $x\in X$, $x's\wedge \neg x's = 0$ is the bottom, and $x's\vee \neg x's =1$ is the top. 

\indent Similarly, we can define the free distributive lattice $L(P', X)$ for any subset 
$P'\subseteq P$. In particular, we have the free distributive lattice 
$L(\{p\}, X)$ for any $p\in P$.

\indent Since $(L(P,X), \wedge, \vee, ,\leq, 0, 1)$ is a finite distributive lattice, every element can be represented as $M_1\vee M2\vee ... \vee M_k$, where each $M_i$ is a meet of elements from $L_0$. 

%\begin{definition}\label{definition201}
%Let $(P, \leq)$ be a partially ordered set and let $X$ be a set. For any element 
%$p\in (L(P, X),\wedge, \vee, \leq, 0, 1)$, if 
%$p=M_1\vee M2\vee ... \vee M_k$ with each $M_i=\bigwedge S_i$ for some 
%$S_i\subseteq L_0$, then the family $\{S_i: 1\leq i\leq k\}$ is called the support 
%family of $p$, denoted by $supp(p)$. 
%\end{definition}

\begin{definition}\label{definition1}
 (Typed Topological Space)
Let $(X, \sT)$ be a topological space, and let $(P, \leq)$ be a partially ordered set.
A function $\sigma:\sT\rightarrow (L(P, X), \wedge, \vee, \leq, 0, 1)$ is called a type mapping whenever the following conditions hold: 
(1) $\sigma(U)=0$ if and only if $U=\emptyset$; 
(2) $\sigma(U)\neq 1$; and 
(3) for any two open sets $U, V\in\sT$,  the condition $U\subseteq V$ 
implies $\sigma(U)\leq \sigma(V)$.

The 5-tuple $(X, \sT, P, \leq, \sigma)$ is called a typed topological space, 
and each non-empty open set $U$ is called $\sigma(U)$-typed open set, which is
 denoted by $\sigma(U)\Vdash U$. The denotation $\sigma(U)\Vdash x$ is used when $x\in U$.
\end{definition}

\vspace{0.2in}

\indent Certainly, both $(\sT, \cap, \cup, \subseteq)$ and $(L(P, X),\wedge, \vee, \leq)$ are finite lattices. In general, the type mapping is not a lattice morphism which preserves meet and join. 

%\indent Since in a finite partial order set, every filer is principle, i.e., any filter will have the %smallest element, The type-mapping in above definition \cite{definition1} can be rephrased as a %function $\sigma:\sT\rightarrow L(P, X)$ such that for any two open sets $U,V\in\sT$, the condition %$U\subseteq V$ implies $\sigma(U)\leq \sigma (V)$. We use the definition in terms of  filters to %indicate that an open set can be characterized a collection of types, which forms a filter.

%\indent We may also define the type mapping $\sigma$ as in the following: $\sigma(U)\subseteq L_0$, %since any subset of $L_0$ can generate a filter base. To specify a subset of types for a %given open set, we want to emphasize that only those types generated by that subset of types are %direct related to the open set. 

\begin{proposition}\label{proposition22}
Let $(X, \sT, P, \leq, \sigma)$ be a typed topological space. If $U, V\in\sT$, then $\sigma(U\cap V)\leq \sigma(U)\wedge\sigma(V)$, and $\sigma(U)\vee \sigma(V)\leq \sigma(U\cup V)$. $\Box$
\end{proposition}

\indent The concept of typed topology can be applied to situations that involves topologies on finite sets. To the memory of my academic advisor and academic father, W.W. Comfort, we first introduce the following example.

\begin{example}\label{example21}
{\normalfont
(Mathematics Genealogy)  As of December 26, 2017, the math genealogy project collects 221844 records of mathematicians in human history. Let $X$ be the set of all those mathematicians listed there. We will define a typed topology on $X$. Let $P$ be the set of two elements 
$\{ancestors, descendants\}$. For each mathematician 
$x$, we have $x's$ 
and $\neg x's$. So $L_0=\{ancestors, descendants\}\cup\{x's, \neg x's: x\in X\}$. 
\\
\indent According to the semantics of $x's$, ancestors and descendants, we can define open sets on $X$. For instance, the author's name is W. Hu. The open set $U_{W. Hu's, ancestors}$ associated with the types of W. Hu's and ancestors is defined as 
\\
\indent $U_{W.Hu's, ancestors} = 
\{$
W.W.Comfort(1), 
E. Hewitt(2), 
M.H.Stone(3), 
George, D. Birkhoff(4), ...
%E.H.Moore(5), 
%H.A. Newton(6), 
%Michel Chasles(7), 
%Simeon D. Poisson(8), 
%Joseph L. Lagrange(9),
%Pierre-Simon Laplace(9), 
%Leonhard Euler(10), 
%Jean Le Rond, d'Alembert(10), 
%Johann Bernoulli(11), 
%Jacob Bernoulli(12), 
%Nikolaus Eglinger(12), 
%Nicolas Malebranche(13),
%Gottfried Wilhelm Leibniz(14),
%Emmanuel Stupanus(14)
%    Johann Caspar Bauhin(14),
%    Jakob Thomasius(15), 
%Erhard Weigel(15), 
%Petrus Ryff(15),
%Emmanuel Stupanus(15),
%    Friedrich Leibniz(16),
%Philipp Müller(16),
%Felix Plater(16), 
%Christoph Meurer(17), 
%Guillaume Rondelet(17)
%Moritz Valentin Steinmetz(18), 
%Johannes Winter von Andernach(18),
%Georg Joachim von Leuchen Rheticus(19)
%Johann Hoffmann(19), 
%Rutger Rescius(19),
%Jacobus (Jacques Dubois) Sylvius(19),
%Johannes Volmar(20)
%Nicolaus (Mikolaj Kopernik) Copernicus(20),
%    Girolamo (Hieronymus Aleander) Aleandro_{20},
%    Jean Tagault_{20},
%     François Dubois_{20},
%%     Bonifazius Erasmi_{21},  ...
$\}$ 
%\\
%\indent 
%Here, for each name, there is a number followed, e.g., W.W.Comfort(1). It is interpreted that
% W.W.Comfort is the first generation ancestor of W. Hu. Similarly, E. Hewitt(2) is interpreted 
% that E. Hewitt is the second generation ancestor of W. Hu.
In fact, for any direct student $x$ of W.W. Comfort, we have 
$U_{x's, ancestors} = U_{W.Hu's, ancestors}$.
       Similarly, we define $U_{W.W.Comfort's, ancestors}$, and $U_{W.W.Comfort's, ancestors}$
=$U_{W.Hu's, ancestors}$$\setminus$\{W.W.Comfort(1)\}. 
\\
\indent    Similarly, we let $U_{W.Hu's, descendants}$ be the set of all descendants of W. Hu. In this case, it is the empty set, while $U_{W.W.Comfort's, descendants}$ contains 91 names which includes W. Hu. 
    Certainly, $U_{W.W.Comfort's, ancestors}\subseteq U_{W.Hu's, ancestors}$, and $U_{W.Hu's, descendants}\subseteq U_{W.W.Comfort's, descendants}$.

\indent The type mapping $\sigma$ follows from above definition. 
$\sigma(U_{W.Hu's, ancestors})$=
 W.Hu's $\wedge$ ancestors $\wedge\bigwedge$
\{x's: x~is~a~descendant~of~W.W.Comfort\},
and 
$  \sigma(U_{W.Hu's, descendants})=0$ since $U_{W. Hu's, descendants} =\emptyset$ .  
   $\sigma(U_{W.W.Comfort's, ancestors})$
=W.W.Comfort's $\wedge$ ancestors
$\wedge\bigwedge$
\{x's: x~is~a~descendant~of~E.Hewitt\}. 
$\sigma(U_{W.W.Comfort's, descendants})$
=W.W.Comfort's
$\wedge$ descendants
$\wedge\bigwedge$\{x: x~is~an~ancestor~of~W.W.Comfort\}.
 \\
\indent  Pierre-Simon Laplace's advisor is Jean~Le~Rond, d'Alembert, who does not have ancestors in the database. Hence, $U_{Pierre-Simon Laplace's, ancestors}$=\{Jean Le Rond, d'Alembert\} is a singleton that is closed and open. However, when we restrict open sets to those  of types involving ancestors but not descendants, $U_{Pierre-Simon Laplace's, ancestors}$  is not closed anymore. In fact, its closure is the set of all 99910 descendants of Jean Le Rond, d'Alembert. The closure of $U_{Pierre-Simon Laplace's, ancestors}$ using open sets of types involving descendants but not ancestors is itself.
      When a mathematician, e.g., W.Hu has no descendants (students), the singleton $\{W.Hu\}$ is a closed set, since for any other mathematician $x$, $U_{x's, ancestors}$ is disjoint with $\{W. Hu\}$. Whether it is open depends on other students of W.W.Comfort.
      When a mathematician $x$ either has only one advisor and a student, or has an advisor and only one student, the singleton $\{x\}$ is open and close. For instance, $U_{W.Hu's, ancestors}\cap U_{E.Hewitt's, descendants} = \{W.W.Comfort\}$, is an open set. It is closed, since $\{W.W.Comfort\}\cap (U_{W.W.Comfort's, ancestors}\cup U_{W.W.Comfort's, descendants})=\emptyset$. 
%\\  
%\indent One can extends the set $P$ to include more types. For instance, 
%when PhD year after 1800 is added, $U_{W.Hu's, ancestors, PhD~year~after~1800}=\{$
%W.W.Comfort(1), 
%E. Hewitt(2), 
%M.H.Stone(3), 
%George, D. Birkhoff(4), 
%E.H.Moore(5), 
%H.A. Newton(6), 
%Michel Chasles(7), 
%Simeon D. Poisson(8)
%$\}$. 
%The closure of $U_{Pierre-Simon Laplace's, ancestors}$ using open sets of types of 
%ancestors and PhD year after 1800 but not descendants is the set of all  descendants of 
%Jean Le Rond, d'Alembert before year 1800. 
}$\Box$
\end{example}

\begin{example}\label{example22}
{\normalfont
(Community And Neighborhood)  In a community, there are 5 streets. Residents on each street have their neighbors. Let $X$ be the set of all residents in a given community.
Let $P$ be an unordered set that includes 5 street names. Each street name defines an open set, which consists of all residents on that street. In addition, we add two more types into $P$, i.e., "left-neighbor" and "right-neighbor". For a resident on a street $s$, we define $U_{x's, s, left~neighbor}$ to be the set of all residents on the street $s$ that are on the left-hand side of $x$, including $x$, and $U_{x's, s, right~neighbor}$ to be the set of all residents on the street $s$ that are on the right-hand side of $x$ including $x$. The type mapping can be defined naturally. Hence, we have a typed topological space. 

\indent Without considering types, one can show that every singleton is closed and open. When using open sets that involving right-neighbor type but not left-neighbor type, the closure of each singleton $\{x\}$ is the set of all left neighbors of $x$ on the same street. Similarly, when using open sets that involving left-neighbor type but not right-neighbor type, the closure of the singleton $\{x\}$ is the set of all right neighbors of $x$ on the same street.
$\Box$
}
\end{example}

\begin{example}\label{example23}
{\normalfont
(Database)  A relational database uses the so-called relational model to store data. In that model, data is represented as entity and relationships. Each table defines one entity type, of which instances (or records) of that entity type (such as customer, product) are organized as rows. Each instance is further divided into attributes, i.e., columns (fields). Tables are connected by relationships that are defined through some columns such as names shared by both tables. For instance, a table of the entity type customer and a table of the entity type office visit can be related by the same last name and first name. 
\\
\indent Let $X$ be the set of all instances (rows) from all tables in a relational database. Since each table represents an entity type, we will add the table name to the partial order set $P$. Further, each column in a table describes an attribute and each attribute is assigned a domain, which together form a type.
 For instance, an amount column $fee$ may form predicates such as
 "$fee < 1000$" or "$fee\geq 200$". Similarly, a column of product name $product$ may form predicate such as "product likes S*", which means the product name starts with "S". We will include predicates of each attribute as a subset of $P$. The partial order on $P$ is based on the partial order on predicates. 

\indent    When using Structured Query Language (SQL) statements to operate against a database, one can use any of the eight operators, i.e., union, intersection, difference, Cartesian product, selection, projection, join, and division. The selection operator lets us define open sets according to predicates based on attributes. Typically, we use a statement such as 
\\
"SELECT * \\
FROM mytable\\ 
WHERE columnx  $>$ 100"
\\ 
to define a subset of rows in a table. The predicate "columnx > 100" is a type and the result set of rows from the table "mytable" is an open set associated with that type. Basic predicates may use operators such as =, $\neq, >, \geq, <. \leq, IN,$ BETWEEN, LIKE, IS NULL, IS NOT NULL. Compound predicates may connect basic ones by "AND, "OR" and parenthesis.

\indent The definition of open sets automatically defines the type mapping. Hence, we have a typed topological space. Certainly, we can write an SQL statement that selects any given row in a table. Hence, each singleton is an open set and therefore a closed and open set. However, when restricted to a set of predicates, the selected set of rows may not be a singleton. $\Box$
}
\end{example}

\indent For any typed topological space $(X, \sT, P, \leq, \sigma)$, 
since for any non-empty open set $U\in\sT$, $\sigma(U)\neq 0$, we have the following  proposition.

\begin{proposition}
Let $(X, \sT, P, \leq, \sigma)$ be a typed topological space. 
For any $x\in U\in\sT$, if $p, q\in L(P, X)$ are incompatible, then
we cannot have both $p\Vdash x$ and $q\Vdash x$ at the same time. In particular, 
we cannot have both $p\Vdash x$ and 
$\neg p\Vdash x$ at the same time. 
\end{proposition}

{\sl Proof: } Let everything be as above. Assume for a contradiction, there are open sets $W,V\in\sT$ such that $\sigma(W)=p$,  $\sigma(V)=q$ and $x\in W\cap V$. Then by Proposition \ref{proposition22}, we have $\sigma(W\cap V)\leq \sigma(W)\wedge \sigma(V)=p\wedge q$. Since $p$ and $q$ are incompatible, we have $p\wedge q=0$. Hence $\sigma(W\cap V)=0$ and $W\cap V=\emptyset$, which is however a contradiction with $x\in W\cap V$. 
$\Box$

\indent In a finite lattice $L$, any filter $F$ has the smallest element $\bigwedge F$ satisfying $\bigwedge F\leq p$ for any $p\in F$. In general, for any other element $q\in L$, we may not have either $q\in F$ or $\neg q\in F$ (or equivalently, either $\bigwedge F\leq q$ or $\bigwedge F\leq\neg q$). 
When $p$ is join-prime (i.e., $p\leq q\vee r\rightarrow (p\leq q)\vee (p\leq r)$), we have for any $q\in L$, either $p\leq q$ or $p\leq\neg q$ holds. As a known fact, in a distributive lattice, join-prime is equivalent to join-irreducible, which is further equivalent to being atomic.

\begin{lemma}\label{lemma210}
Let $(L, \wedge, \vee, \leq, 0, 1)$ be a finite distributive lattice. For any join-irreducible element $p\in L$, the family $\{q\in L: p\leq q\}$ is an ultrafilter containing p. For any ultrafilter $F$ in $L$, $\bigwedge F$ is a join-irreducible element. 
\end{lemma}

{\sl Proof: } When $p$ is join-irreducible, or equivalently join-prime, the family $F=\{q: p\leq q\}$ is a filter, since $p\leq q\wedge r$ for any two elements $q,r\in F$. It is also an ultrafilter, since for any $q\in L$, $p\leq q\vee\neg q=1$, which implies either $q\in F$ or $\neg q\in F$. 

\indent When $F$ is an ultrafilter in $L$, for any $q\in L$, we have either $q\in F$ or $\neg q\in F$, which implies $\bigwedge F\leq q$ or $\bigwedge F\leq\neg q$ respectively. Hence $\bigwedge F$ is a join-irreducible (or join-prime) element. $\Box$

\indent When a filter $F$ is not an ultrafilter, we may have element $p\in P$ satisfying neither $p\in F$ nor $\neg p\in F$. The following proposition is easy to verify.

\begin{lemma}
Let $(X, \sT, P, \leq, \sigma)$ be a typed topological space. For any $U\in\sT$, if $\sigma(U)$ is join-irreducible, then for any  $p\in L(P,X)$ either $\sigma(U)\leq p$ or $\sigma(U)\leq \neg p$. 
    Furthermore, the pair $\{p, \neg p\}$ can be replaced by any maximum antichain in $(L(P, X), \leq)$.
$\Box$
\end{lemma}

\section{Typed neighborhood systems}

\indent\indent As we mentioned above, most singletons in a finite typed topological space are closed and open. However, when limiting the open sets to certain types, they may not be that case anymore. In applications such as a database, we usually want to find a set of points (records) by a minimum set of types. Let us investigate the neighborhood systems of a point in $(X, \sT)$.

\indent Recall, in a finite distributive lattice $(L, \wedge, \vee)$ with bottom and top, an element $p\in L$ is called $join$-$irreducible$ if $p\neq 0$ and whenever $p= q\vee r$ one must have either $p=q$ or $p=r$, and $meet$-$irreducible$ if $p\neq 1$ and whenever $p= q\wedge r$ one must have either $p=q$ or $p=r$. As a known fact, every element in such type of lattice is a meet of meet-irreducible elements or a join of join-irreducible elements. 

\indent Let $(X, \sT, P, \leq, \sigma)$ be a typed topological space. Then, $(\sT, \cap, \cup)$ is a finite distributive lattice. The family of all join-irreducible open subsets of $X$ is a base of $\sT$. 

\begin{definition}\label{definition250}
 Let $(X, \sT, P, \leq, \sigma)$ be a typed topological space. An open set $U\in\sT$ is called $p$-$join$-$irreducible$ for $0\neq p\in L(P,X)$ if $U\neq \emptyset$, $p\leq\sigma(U)$ and whenever $U=W\cup V$ for some $U, V\in\sT$ satisfying $p\leq \sigma(W)$ and $p\leq \sigma(V)$, one must have either $U=W$ or $U=V$. Similarly, we can define  the concept of being $p$-$meet$-$irreducible$. 
\end{definition}

\indent Certainly, if an open set $U$ satisfying $p\leq\sigma(U)$ is $p$-join-irreducible then it is also $q$-join-irreducible for any $q\in L(P, X)$ with $p\leq q$. Furthermore, if an open set $U$ is
join-irreducible, then it is $p$-join-irreducible for all $p\in L(P, X)$. A $p$-join irreducible open set may not be join-irreducible. 

\begin{example}\label{example31}
{\normalfont
In the Example \ref{example21} of Mathematics Genealogy,  $\sigma(U_{W.Hu's, ancestors})$ is
equal to $W.Hu's \wedge ancestors\wedge\bigwedge$\{x's: x~is~a~descendant~of~W.W.Comfort\}.
 That open set is the only one whose type is 
$W.Hu's \wedge ancestors\wedge\bigwedge$\{x's: x~is~a~descendant~of~W.W.Comfort\}.
Let $p_0=\sigma(U_{W.Hu's, ancestors})$.
If $q\in L(P, X)$ and $p_0\lneq q$, then either $q$ is a meet of some elements from 
\{W.Hu's, ancestors\}$\cup$\{x's: x~is~a~descendant~of~W.W.Comfort\} or 
$q$ is a join of $q'$ 
and some $r$ not in \{W.Hu's, ancestors\}$\cup$\{x's: x~is~a~descendant~of~W.W.Comfort\}, where $q'$ is a meet of some element from 
\{W.Hu's, ancestors\} $\cup$\{x's: x~is~a~descendant~of~W.W.Comfort\}. 
According to Proposition \ref{proposition22}, there is no open set $W$ satisfying $W\subseteq U$ and $\sigma(W)=q$.  Hence, $U_{W.Hu's, ancestors}$ is $p_0$-join-irreducible. 
Similarly, $U_{E.Hewitt's, descendants}$ is $p_1$-join-irreducible, 
where $p_1=\sigma(U_{E.Hewitt's, descendants})$. 
\\
\indent However, if we set $p=p_0\wedge p_1\wedge\sigma(U_{W.W.Comfort's, ancestors})$, 
then $U_{W.Hu's, ancestors}$ is not $p$-join-irreducible, since 
$U_{W.Hu's, ancestors}$ is the disjoint union of two open sets \{W.W.Comfort\} and 
$U_{W.W.Comfort's, ancestors}$, where both  
$p\leq \sigma(\{W.W.Comfort\})$ and 
$p\leq\sigma(U_{W.W.Comfort's, ancestors})$ are true.

\indent On the other hand, the singleton $\{W.W.Comfort\}$ is an open set. Certainly, it is $p$-join-irreducible for any $p\in L(P, X)$. $\Box$
}
\end{example}

\indent The type mapping in a typed topological space does not guarantee that
$\sigma(U)\lneq\sigma(V)$
when $U$ is a proper subset of $V$. 

\begin{definition}\label{definition33}
Let $(X, \sT, P, \leq, \sigma)$ be a typed topological space. When 
$\sigma(U) \lneq \sigma(V)$ holds for any two non-empty 
open sets $U$ and $V$ with $U\subsetneq V$, the space is called strictly typed. 
\end{definition}

\indent The spaces in Example \ref{example21}, \ref{example22} and \ref{example23} are strictly typed topological spaces. In the rest of this article, we will assume all spaces are strictly typed.

\begin{lemma}\label{lemma331}
Let $(X, \sT, P,\leq, \sigma)$ be a strictly   typed topological space. 
For any non-empty open set $U\in\sT$, the set $U$ is $\sigma(U)$-join-irreducible. 
\end{lemma}

{\sl Proof: } For any two open sets $V, W\in\sT$ satisfying $U=W\cup V$, if both $W$ and $V$ are proper subsets of $U$, then by the assumption of being strictly typed, we have both 
 $\sigma(V) \lneq \sigma(U)$ and $\sigma(W) \lneq \sigma(U)$.
Hence $U$ cannot be the union of two non-empty proper open sets $W,V$ satisfying 
$\sigma(U)\leq\sigma(W)$ and $\sigma(U)\leq\sigma(V)$, i.e., $U$ 
is $\sigma(U)$-join-irreducible.
 $\Box$.

\begin{definition}\label{definition331} 
Let $(X, \sT, P, \leq, \sigma)$ be a strictly typed topological space. For any $p\in L(P, X)$ set  $\sT_{\geq p}(X)=\{U\in\sT: p\leq \sigma(U)\}$, and 
$J_{\geq p}(X)=\{U\in\sT_{\geq p}(X): U$ is $p$-join-irreducible$\}$. 
      Furthermore, for any $x\in X$, set $\sT_{\geq p}(x)=\{U\in\sT_{\geq p}(X): x\in U\}$ and
$J_{\geq p}(x)=\{U\in J_{\geq p}(X): x\in U\}$.
\end{definition}

\begin{theorem}\label{theorem32}
Let $(X, \sT, P, \leq, \sigma)$ be a strictly typed topological space. For any $p\in L(P, X)$, every open set $U\in \sT_{\geq p}(X)$ is a join of some elements from $J_{\geq p}(X)$.
\end{theorem}

{\sl Proof: } Since $L(P, X)$ and $\sT$ are finite, we can prove it by induction on open sets. 
We first consider three base cases. Case 1: for any open set $U\in \sT_{\geq p}(X)$ satisfying 
$\sigma(U)=p$, by Lemma \ref{lemma331}, $U$ is $p$-join-irreducible. Hence $U\in J_{\geq p}(X)$. 
Case 2: for any $U\in\sT_{\geq p}(X)$, if there exists no open set $V\in\sT_{\geq p}(X)$ 
satisfying $V\subsetneq U$, then $U$ is  $p$-join-irreducible, i.e., $U\in J_{\geq p}(X)$. 
Case 3: for any $U\in\sT_{\geq p}(X)$, if there do not exist two distinct open sets 
$W, V\in\sT_{\geq p}(X)$ satisfying $U=W\cup V$ and $(W\neq U)\wedge (V\neq U)$, then 
$U$ is $p$-join-irreducible, i.e., $U\in J_{\geq p}(X)$. 
In all three base cases, $U$ is an element in $J_{\geq p}(X)$.

\indent At the induction step, let us assume that for a given open set $U\in \sT_{\geq p}(X)$ satisfying $p\leq \sigma(U)$, any $V\in\sT_{\geq p}(X)$ with $V\subsetneq U$ and
 $p\leq\sigma(V)\lneq\sigma(U)$ is a join of $p$-join-irreducible open sets. We will show 
that $U$ is also a join of $p$-join-irreducible open sets. If there are no $W, W'\in\sT_{\geq p}(X)$
 satisfying $U=W\cup W'$ and 
$p\leq \sigma(W)\wedge \sigma(W')$, then $U\in J_{\geq p}(X)$. Otherwise, let $U=W\cup W'$ with 
 $p\leq \sigma(W)\wedge \sigma(W')$. By assumption, both $W$ and $W'$ are joins of 
$p$-join-irreducible open sets. Therefore, $U$ is also the join of $p$-join-irreducible open sets. 
$\Box$

\indent Neighborhoods of a point can be organized by chains of types, a way to sandwich open sets. By sandwiching open sets' types, we are able to limit those types from $P$ in the operation.

\begin{definition}\label{definition333001}
Let $(X, \sT, P, \leq, \sigma)$ be a strictly typed topological space. 
For $k>1$, let $c= "p_0\leq p_1\leq ... \leq p_{k-1}"$ be a chain of types in $L(P,X)$. For any point $x\in X$, the  $c-$neighborhood system is defined as  the family  $\sT_c(x)=\{U\in\sT: (x\in U)\wedge \exists i(0\leq i<k-1) p_i\leq \sigma(U)\leq p_{i+1}\}$. 
     Furthermore, we define 
$\sJ_c(x)=\{U\in\sT_c(x): U\in J_{\geq p_i}(x)~for~some ~0\leq i<k-1\}$.
We also set  $\sT_c(X)=\bigcup\{\sT_c(x): x\in X\}$, 
and $\sJ_c(X)=\bigcup\{\sJ_c(x): x\in X\}$.
\end{definition}

\begin{theorem}\label{theorem333}
Let everything be as in Definition \ref{definition333001}. 
 Then $\sJ_c(x)$ is a neighborhood base of the $c$-neighborhood system $\sT_c(x)$. 
\end{theorem}

{\sl Proof: } For any $U\in\sT_c(x)$, we need to show that there exists $V\in\sJ_c(x)$ such that $x\in V\subseteq U$. For this $U$, there exists $i$, $0\leq i<k-1$, such that $p_i\leq\sigma(U)\leq p_{i+1}$. Hence $U\in \sT_{\geq p_i}(X)$. By Theorem \ref{theorem32}, $U$ is the join of elements from $J_{\geq p_i}(X)$, say $U=V_1\cup ...\cup V_m$ with $V_k\in J_{\geq p_i}(X)$ for each $k\leq m$. Since $x\in U$, we have $x\in V_j$ for some $j\leq m$. Hence $V_j\in J_{\geq p_i}(x)$.
To complete the proof, we will show that $V_j\in\sT_c(x)$, for which it suffices to show that 
$V_k\in\sT_c(x)$ for all $k\leq m$. Since $p_i\leq \sigma(U)\leq p_{i+1}$, we have $\sigma(V_k)\leq p_{i+1}$ for all $k\leq m$. On the other hand, we have $p_i\leq \sigma(V_k)$ since $V_k\in J_{\geq p_i}(X)$. 
Hence $p_i\leq\sigma(V_k)\leq p_{i+1}$, i.e., $V_k\in\sT_c(x)$.
$\Box$

\indent While $\sT_{\geq p}(X)$ is the family of open sets whose types are above $p$, i.e., the family of open sets of upper types, we can also define the family of open sets whose types are below $p$. However, that family allows any kind of intersections with types including $p$, which makes it less interesting since more closed and open singletons are included. The following definition provides a modified version of lower types than $p$. 

\begin{definition}\label{definition3331}
Let $(X, \sT, P, \leq, \sigma)$ be a strictly typed topological space. Let $p\in P$ be a type. A chain of types $c= "p_0\leq p_1\leq ... \leq p_{k-1}"$ is called a $p$-chain, if each $p_i\in L(\{p\}, X)$ and $p_{k-1}=p$. For $x\in X$, we define the $p-chain$-neighborhood system as $\sT_{p-chain}(x)=\bigcup\{\sT_c(x):$ c is a p-chain $\}$. 
%and $\sJ_{p-chain}(x)=\bigcup\{\sJ_c(x): $ c is a p-chain $\}$. 
Further, we set $\sT_{p-chain}(X)=\bigcup\{\sT_{p-chain}(x): x\in X\}$ 
%and $\sJ_{p-chain}(X)=\bigcup\{\sJ_{p-chain}(x): x\in X\}$.
\end{definition}

\indent In the Example \ref{example21}, let $p$ be the type $ancestors$. Then the family $\sT_{p-chain}(X)$ is family of all open sets of the form $U_{x's, ancestors}$.

\begin{proposition}\label{proposition3331}
Let everything be as in Definition \ref{definition3331}. The $p-chain$-neighborhood system $\sT_{p-chain}(x)$ is equal to the
 family $\{U\in\sT: x\in U\wedge \sigma(U)\in L(\{p\}, X) \wedge \sigma(U)\leq p \}$. 
\end{proposition}

{\sl Proof: } For any $q\in L(\{p\}, X)$ and $q\leq p$, we can form a chain $q\leq p$, which is a $p$-chain. If an open set 
$U\in \{U\in\sT: x\in U\wedge \sigma(U)\in L(\{p\}, X) \wedge \sigma(U)\leq p \}$, then 
$U\in \sT_c(x)$, where $c= "\sigma(U)\leq p"$. Hence $U\in\sT_{p-chain}(x)$. The other direction of inclusion is obvious from the definition of $p$-chain. $\Box$

\begin{theorem}\label{theorem33441}
Let everything be as in Definition \ref{definition3331}. For any $x\in X$, every open set $U\in\sT_{p-chain}(x)$ is an element in $\sJ_c(x)$ for some $p$-chain $c$. Hence,
$\sT_{p-chain}(x)=\bigcup\{\sJ_c(x): $ c is a p-chain $\}$. 

\end{theorem}

{\sl Proof: } Set $p_0=\sigma(U)$, and let $c= "p_0\leq p"$. Then, by Lemma \ref{lemma331}, $U$ is $p_0$-join-irreducible. Hence $U\in\sJ_c(x)$. $\Box$

%\indent The following corollary is straightforward from above Proposition 
%\ref{proposition3331} and Theorem \ref{theorem333}.

%\begin{corollary}\label{corollary33441}
%Let everything be as in above definition. For any $x\in X$ and 
%$p\in L(P,X)$, $\sJ_{p-chain}(x)$ is a neighborhood base of the $p-chain$-neighborhood 
%system $\sT_{p-chain}(x)$. $\Box$
%\end{corollary}

\indent For a finite distributive lattice $L$, it is a known fact that it can be partitioned into $w(L)$-many chains, where $w(L)$ is the largest size of maximum antichain. Hence, there exists a family of $w(L)$-many chains $\{c_i: i<w(L)\}$ that covers $L$, i.e., $L=\bigcup\{c_i: i<w(L)\}$. 

\begin{theorem}\label{theorem334}
Let $(X, \sT, P, \leq, \sigma)$ be a strictly typed topological space. Let $\{c_i: i<w(L(P,X))\}$ be a family of chains that covers $L(P,X)$. Then, for any $x\in X$, $\bigcup\{\sT_{c_i}(x): i<w(L(P,X))\}$
is the neighborhood system of $x$, and $\bigcup\{\sJ_{c_i}(x): i<w(L(P,X))\}$ is a neighborhood base.
\end{theorem}

{\sl Proof: } Since $\{c_i: i<w(L(P,X))\}$ covers $L(P,X)$, for any $U\in\sT$ satisfying $x\in U$, there exists $i<w(L(P,X))$ such that $\sigma(U)\in c_i$. Then $U\in\sT_{c_i}(x)$. Hence 
$\bigcup\{\sT_{c_i}(x): i<w(L(P,X))\}$ is the neighborhood system of $x$. By Theorem \ref{theorem333},  $\sJ_{c_i}(x)$ is a neighborhood base of $\sT_{c_i}$. Hence, there exists $V\in\sJ_{c_i}(x)$ satisfying $x\in V\subseteq U$. Hence $\bigcup\{\sJ_{c_i}(x): i<w(L(P,X))\}$ is a neighborhood base of $x$.
$\Box$

\section{Typed Closure}

\indent\indent As we discussed before, closure in a typed topological space will be more meaningful when it is restricted to types. 

\begin{definition}\label{definition341}
Let $(X, \sT, P, \leq, \sigma)$ be a strictly typed topological space. For any subset $A\subseteq X$
 and any chain $c$ of types, the $c$-closure of $A$, denoted $\overline{A}^c$, is the set 
$\{x\in X: (\forall U\in\sJ_c(x)) U\cap A\neq\emptyset\}$.
\end{definition}

\indent Since $\sJ_c(x)$ is finite, the following lemma is straightforward.

\begin{lemma}\label{lemma341}
Let everything be as in Definition \ref{definition341}. For any $x\in X$ with $\sJ_c(x)\neq\emptyset$, $x\in\overline{A}^c$ if and only if $\bigcap\sJ_c(x)\cap A\neq\emptyset$. $\Box$
\end{lemma}

\indent If $x\in\overline{A}^c$ and $\sJ_c(x)\neq\emptyset$, then $\bigcap\sJ_c(x)\cap A\neq\emptyset$. 
Let $y\in \bigcap\sJ_c(x)\cap A$. Then $\sJ_c(x)\subseteq\sJ_c(y)$. We have following theorem.

\begin{theorem}\label{theorem3401}
Let everything be as in Definition \ref{definition341}. For any two distinct points $x,y\in X$ satisfying $\sJ_c(X)\neq\emptyset$ and $\sJ_c(y)\neq\emptyset$, the following conditions are equivalent
\begin{enumerate}
\item $\sJ_c(x)\subseteq\sJ_c(y)$; 
\item $x\in\overline{A}^c$ holds for all subsets $A\subseteq X$ with $y\in A$; and
\item $x\in\overline{\{y\}}^c$.

\end{enumerate}
\end{theorem}

{\sl Proof: } (1)$\rightarrow $(2). When $\sJ_c(x)\subseteq\sJ_c(y)$, we have $y\in\bigcap\sJ_c(x)$, which implies $\bigcap\sJ_c(x)\cap A\neq\emptyset$. Hence by Lemma \ref{lemma341}, (2) holds.
\\
\indent (2)$\rightarrow$(3). Trivial. 
\\
\indent (3)$\rightarrow$(1). By Lemma \ref{lemma341} again, $\bigcap\sJ_c(x)\cap\{y\}\neq\emptyset$, which means $y\in\bigcap\sJ_c(x)$. Hence (1) holds. $\Box$

\begin{definition}\label{definition3401} 
Let everything be as in Definition \ref{definition341}. For any chain $c$ of types, set $E_c(X)=\{x\in X: \sJ_c(x)=\emptyset\}$. 
\end{definition}

\indent We have following result.

\begin{lemma}\label{lemma3411}
$E_c(X)=X\setminus\bigcup\{\bigcup\sJ_c(x): x\notin E_c(X)\}$, and $E_c(X)$ is a $c$-closed subset of $(X, \sT)$. 
\end{lemma}

{\sl Proof: } By definition, for any $x\notin E_c(X)$, and any $U\in\sJ_c(x)$, we have $U\cap E_c(X)=\emptyset$. Hence $\bigcup\sJ_c(x)\cap E_c(X)=\emptyset$. 
$\Box$

\begin{definition}\label{definition342}
Let everything be as above. For any subset $Y\subseteq X$, a set $D\subseteq Y$ is called $c$-dense in $Y$ if for any point $x\in Y$ either $x\in D$ whenever $x\in E_c(X)$ or $U\cap D\neq\emptyset$ holds for all $U\in\sJ_c(x)$. 
\end{definition}

\begin{lemma}\label{lemma3412}
Let everything be as above. A set $D$ is $c$-dense in $Y$ if $E_c(X)\cap Y\subseteq D$ and for any point $x\in Y\setminus E_c(X)$, there exists $y\in D$ such that $\sJ_c(x)\subseteq\sJ_c(y)$.
\end{lemma}

{\sl Proof: } Certainly, $E_c(X)\cap Y\subseteq D$ must be true for $D$ to be $c$-dense in $A$. For $x\in Y\setminus E_c(X)$, if $Y\cap \bigcap\sJ_c(x)=\{x\}$, then $x\in D$ must be true; otherwise 
$\bigcap\sJ_c(x)\cap D=\emptyset$, which contradicts with the assumption that $D$ is $c$-dense in $Y$.

When $Y\cap \bigcap\sJ_c(x)$ contains elements other than $x$, we have either $x\in D$ 
or there exists $y\in \bigcap\sJ_c(x)\cap D$. In the first case, our choice of $y$ will be 
$x$ itself, which satisfying $\sJ_c(x)\subseteq\sJ_c(y)$. In the second case, 
we also have $\sJ_c(x)\subseteq\sJ_c(y)$, since $y\in\bigcap\sJ_c(x)$.  $\Box$

\indent Set $\sJ_c=\{\sJ_c(x): x\notin E_c(x)\}$. 
We can order the family $\sJ_c$  by $\subseteq$, which becomes a partially ordered set. 

It may happen that $\sJ_c(x) = \sJ_c(y)$ for $x\neq y\in X$. We can define an equivalent relation $\equiv_c$ on $X\setminus E_c(X)$ as $x\equiv_c y$ whenever $\sJ_c(x) = \sJ_c(y)$. 

\begin{lemma}\label{lemma3413}
Let everything be as above. For any $x\in X\setminus E_c(X)$, if $\sJ_c(x)$ is a maximal element in the partially ordered set $(\sJ_c, \subseteq)$, then for any $A\subseteq X$, $x\in\overline{A}^c$ if and only if 
there exists a $y\equiv_c x$ satisfying $y\in A$.
\end{lemma}

{\sl Proof: } By Lemma \ref{lemma3412}, $x\in\overline{A}^c$ if and only if there exists $y\in A$
satisfying $\sJ_c(x)\subseteq\sJ_c(y)$. Since $\sJ_c(x)$ is a maximal element in the partially ordered set $(\sJ_c, \subseteq)$, we have $\sJ_c(x) = \sJ_c(y)$. Hence $x\equiv_c y$ holds. 
\\
\indent When there exists a $y\equiv_c x$ satisfying $y\in A$, we have $\sJ_c(x)=\sJ_c(y)$ and 
$y\in\bigcap\sJ_c(x)\cap A\neq\emptyset$, which implies $x\in\overline{A}^c$.
$\Box$

\begin{definition}\label{definition3413}
The $c$-density of $(X, \sT)$, denoted $d_c(X, \sT)$, is the smallest size of $c$-dense subset of $(X,\sT)$.
\end{definition}

\begin{theorem}\label{theorem3421}
Let everything be as above. Then $d_c(X, \sT)$ is equal to the sum of $|E_c(X)|$ 
and the number of maximal elements in  $(\sJ_c, \subseteq)$.
\end{theorem}

{\sl Proof: } Let $D$ be a $c$-dense subset of $(X, \sT)$. Then $E_c(X)\subseteq D$. By Lemma \ref{lemma3412}, for any $x\in X\setminus E_c(X)$ with $\sJ_c(x)$ being a maximal element in $(\sJ_c, \subseteq)$, since $x\in X=\overline{D}^c$, there exists a $y\equiv_c x$ satisfying $y\in D$. Hence $|D|$ is greater than and equal to the sum of $|E_c(X)$ and the number 
of maximal elements in $(\sJ_c, \subseteq)$. 

\indent If there exists $x\in D\setminus E_c(X)$ such that $\sJ_c(x)$ is not a maximal element 
in $(\sJ_c, \subseteq)$, then there exists $y\in D\setminus E_c(X)$ such that $\sJ_c(y)$ is maximal in 
$(\sJ_c, \subseteq)$ and $\sJ_c(x)\subseteq \sJ_c(y)$, which implies that $D\setminus\{x\}$ is also 
$c$-dense in $(X, \sT)$. However, that
contradicts with the definition of $d_c(X, \sT)$ being the smallest size of $c$-dense subsets. Hence, every $x\in D\setminus E_c(X)$ satisfies that $\sJ_c(x)$ is a maximal element in $(\sJ_c, \subseteq)$, and the conclusion of the theorem holds. $\Box$ 

\begin{corollary}\label{corollary3422}
Let everything be as above. The $c$-dense subset of $(X,\sT)$ of size $d_c(X, \sT)$ is unique up to the equivalent relation $\equiv_c$.
\end{corollary}

{\sl Proof: } For each $x$ with $\sJ_c(x)$ being a maximal element in $(\sJ_c, \subseteq)$, we pick  one point in its equivalent classes of $\equiv_c$, and form a set $D'$. Then,  
by Theorem \ref{theorem3421}, the smallest size $c$-dense subset $D$ is the union of $E_c(X)$
and $D'$. $\Box$

\section{Typed connected sets, connections and statistics}

%\begin{definition}\label{definition34301}
%       For $k>1$, let $c=p_0\leq p_1\leq ... \leq p_{k-1}$ be a chain of types 
%in $L(P,X)$. Set  $\sT_c(X)=\bigcup\{\sT_c(x): x\in X\}$, 
%and $\sJ_c(X)=\bigcup\{\sJ_c(x): x\in X\}$.      
%\end{definition}

\begin{definition}\label{definition343}
Let $(X, \sT, P, \leq, \sigma)$ be a strictly typed topological space. For any chain of types $c$, a subset $A\subseteq X$ is called $c$-connected if there do not exist  two   disjoint open sets $U,V\in\sT_c(X)$ satisfying $A = A\cap(U\cup V)$, $A\cap U\neq\emptyset$ and 
$A\cap V\neq\emptyset$.  
\end{definition}

\begin{theorem}\label{theorem3333}
Let everything be as a Definition \ref{definition343}. For any $U\in\sJ_c(X)$, if $U$ is $p_0$-join-irreducible, then $U$ is $c$-connected, where $c= "p_0\leq p_1\leq ...\leq p_{k-1}"$.
\end{theorem}

{\sl Proof: } By definition, if $U$ is $p_0$-join-irreducible, then there do not exist any two disjoint non-empty open sets $W, V\in\sT_{\geq p_0}(X)$ satisfying $U=W\cup V$. 
Since $\sT_c(X)\subseteq\sT_{\geq p_0}(X)$, we conclude that $U$ is $c$-connected. 
$\Box$

\begin{theorem}\label{theorem33331}
Let $c= "p_0\leq p_1\leq...\leq p_{k-1}"$ be a chain of types. Then every open set 
$U\in\sJ_{\geq p_0}(X)$ is $c$-connected.
\end{theorem}

{\sl Proof: } By Definition \ref{definition331}, $U\in\sJ_{\geq p_0}(X)$ if $U$ is $p_0$-join-irreducible. Hence, there do not exist two non-empty open sets $W, V\in\sT$ satisfying $p_0\leq\sigma(W)\wedge\sigma(V)$ 
and $U=W\cup V$, which implies that there do not exist two disjoint non-empty open sets $W, V\in\sT_c(X)$ satisfying $U=U\cap(W\cup V)$, $U\cap W\neq\emptyset$ and $U\cap V\neq\emptyset$. Hence $U$ is $c$-connected. $\Box$

\begin{theorem}\label{theorem33332}
Let $p\in P$ be a type. For any $U\in\sT_{p-chain}(X)$, there exists a $p$-chain $c$ such that $U$ is $c$-connected. 
\end{theorem}

{\sl Proof: } By Theorem \ref{theorem33441}, every $U\in\sT_{p-chain}(X)$ is an element in $\sJ_c(X)$ for the $p$-chain $c="\sigma(U)\leq p"$. Hence $U$ is $\sigma(U)$-join-irreducible, i.e., 
$U\in\sJ_{\geq\sigma(U)}(X)$. By Theorem \ref{theorem33331}, $U$ is $c$-connected. $\Box$

%\begin{definition}\label{definition344}
%Let $c_1, c_2, ..., c_n$ be a finite sequence of chains of types. An element $x\in X$ is called %$\{c_i\}$-networked to an element $y\in X$, if there exists elements $x_1, x_2, ..., x_{n+1}\in X$ 
% and a family of subsets $\{A_i: 1<=i<=n\}$ of $X$ such that (1) $x=x_1, y=x_{n+1}$, (2) $A_i$ is %$c_i$-connected and $x_i, x_{i+1}\in A_i$ hold for all $i<=n$. 
%An element $x\in X$ is called a $\{c_i\}$-suspect for an element $y\in X$, if there exists a family %of subsets $\{A_i: 1<=i<=n\}$ such that $A_i$ is $c_i$-connected and $x, y\in A_i$ hold for all $i$.
%\end{definition}

\begin{definition}\label{definition345}
Let everything be as above. Two elements $x, y\in X$ are called $c$-connected if there exists a $c$-connected subset $A$ satisfying $x,y\in A$. The set $A$ is called a connection between $x$ and $y$.
\end{definition}

\begin{example}\label{example345}
{\normalfont
(Community and Neighborhood Revisit)  In Example \ref{example22}, we defined a typed topological space of community and neighborhood. In this example, we revisit it. In addition to the types such as street names, and left-neighbor and right-neighbor, we add the type "classmate". For a resident $x$ on a street $s$, the open set $U_{x's, s, left-neighbor}$ and  $U_{x's, s, right-neighbor}$ are the same. For the type of classmate, we define $U_{x's, classmate}$ to be the set of all residents in the community, including $x$, who are classmates of $x$ in any grades. 

\indent If two elements $x$ and $y$ are classmates, then let $c$ be the chain $"x's\wedge classmate \leq classmate"$. Then, $O=U_{x's, classmate}$ is a $c$-connected subset, and  $x$ and $y$ are $c$-connected by $O=U_{x's, classmate}$. $\Box$

% by the three types. For instance, a person $x_2$ on the same street $s$ of $x$ is a left neighbor %of $x$. The individual $x_2$ has a classmate $x_3$ who lives at street $t$ and $x_3$ is a right %neighbor of $y$. In this case, $c_1$ is the chain of types 
%$x's\wedge s\wedge left-neighbor\leq x's\wedge s\wedge left-neighbor$, 
%$c_2$ is $x_2's\wedge classmate\leq x_2's\wedge classmate$, and 
%$c_3$ is $y's\wedge t\wedge right-neighbor\leq y's\wedge t\wedge right-neighbor$. The subsets are 
%$A_1 =U_{x's,s, left-neighbor}$, 
%$A_2=U_{x_2's, classmate}$, 
%$A_3=U_{y's, t, right-neighbor}$. 
%Certainly, each $A_i$ is $c_i$-connected, since that is the only open set of the type. Hence $x$ is %$\{c_i\}$-networked to $y$. 

}
\end{example}

\indent Since our spaces are finite, statistics can play a key role in studying individual points. For instance, we can measure the neighborhood system of each individual point. We can also measure the relationship between two points. 

\indent Let $c$ be a chain of types, Then, the set $\{|A|: A\subset X$ is a $c$-connected subset$\}$ is a finite set of non-negative integers. According to basic statistics, it has mean and standard deviation. In Example \ref{example345}, for any point $x$ and the chain $c = "x's\wedge classmate\leq classmate"$, there are few (actually 1) sets that are $c$-connected. That may not be enough for statistical calculation. We can include more sets in to the calculation.

\begin{definition}\label{definition346}
Let $(X, \sT, P, \leq, \sigma)$ be a typed topological space. Let $p\in P$ be a type. Then, the sample mean $\overline{x}_p$ and sample standard deviation $s_p$ of $p$ is defined as the mean and standard deviation of the set $\{|A|: A\in\sT_{p-chain}(X)\}$. Furthermore, for any member $A$ in $\sT_{p-chain}$, its $z_p-score$ (or $p$-standard score) is the z-score of $|A|$ with respect to $\overline{x}_p$ and $s_p$.
\end{definition}

\begin{example}
{\normalfont
 In Example \ref{example345}, let $p$ be the type "classmate". The family $\sT_{p-chain}(X)$ 
is the family that contains $c$-connected subsets with $c= "x's\wedge classmate\leq classmate"$ 
for any $x\in X$. Hence, it is the family that contains subsets of the form $\{y\in X$: y is one of x's classmates in the
 community$\}$ for any resident $x$. In that case, $\overline{x}_p$ represents the sample mean of number
 of classmates residents has in the community and $s_p$ represents the corresponding standard
 deviation. 
    Furthermore, for any resident $x$, the corresponding $z_p$-score of $U_{x's, classmate}$ describes
 how many classmates of $x$ in the community, comparing to others in the community. $\Box$ 
}
\end{example}

\indent Individual points in a typed topological space can be evaluated by related statistics. 

\begin{definition}\label{501}
Let $(X, \sT, P, \leq, \sigma)$ be a strictly typed topological space. Given a type $p\in P$,
 for any $x\in X$, set 
$L_x=\{p\in L(\{p\}, X): p\Vdash x\}$. Then the sample mean $\overline{x}_x$ and sample standard deviation 
$s_x$ of the point $x$ are defined as the mean and standard deviation of the set $\{|L_x|: x\in X\}$. Furthermore, the $z_L$-score of $x$ is the $z$-score of $|L_x|$ with respect to $\overline{x}_x$ and $s_x$. 
\end{definition}

\indent Semantically, points with higher $z_L$-scores indicate more active in the space. 
In the Example \ref{example21} of Mathematics Genealogy, mathematicians who have most descendants have much higher $z_L$-scores. 
Similarly, in the Example \ref{example22} of Community and Neighborhood, residents living on two ends of a street have higher $z_L$-scores for an appropriate type. 

\indent The relationship between two points can be measured statistically too. 

\begin{definition}\label{definition55}
Let $(X, \sT, P, \leq, \sigma)$ be a strictly typed topological space. For any two disjoint points 
$x, y$, we define $L_{x, y}=\{p\in L(P, X): p\Vdash x~and ~ p\Vdash y\}$. Then the sample mean $\overline{x}_{x, y}$ and sample standard deviation 
$s_{x, y}$ of the pair $\{x,y\}$ are defined as the mean and standard deviation of the set $\{|L_{x, y}|: x, y\in X\}$. Furthermore, the $z_L$-score of $\{x, y\}$ is the $z$-score of $|L_{x, y}|$ with respect to $\overline{x}_{x, y}$ and $s_{x, y}$.
\end{definition}

\indent Pairs with higher $z_L$-scores indicate more interactions between them. In the Example \ref{example22} of Community and Neighborhood, if we add more types such as "party together", 
"on a team of sports", "sleep over", etc., a pair of two points who has higher $z_L$-score is closer than other pairs. The semantics comes in naturally!

%\section{Lattice morphisms}
% When $\sigma$ is a lattice mrophism, each $\sB_c(x)$ is a topology. 

%\section{Continuous mappings between two typed topological spaces}

\end{document}